\documentclass{amsart}

\theoremstyle{definition}

\theoremstyle{remark}

\numberwithin{equation}{section}

\begin{document}

\title{Cheng Equation: A Revisit Through Symmetry Analysis}


\author{Amlan K Halder}
\address{Department of Mathematics, Pondicherry University, Kalapet, India-605014}
\curraddr{}
\email{amlan.haldar@yahoo.com}
\thanks{AH expresses grateful thanks to UGC (India), NFSC, Award No. F1-17.1/201718/RGNF-2017-18-SC-ORI-39488 for financial support and Late Prof. K.M.Tamizhmani for his encouragement and support.}

\author{R Sinuvasan}
\address{Department of Mathematics, Shanmugha Arts Science Technology and Research Academy, Thanjavur, India 613401}
\curraddr{}
\email{rsinuvasan@gmail.com}
\thanks{}

\author{A Paliathanasis}
\address{Instituto de Ciencias F\'{\i}sicas y Matem\'{a}ticas, Universidad Austral de
Chile, Valdivia, Chile}
\address{Institute of Systems Science, Durban University of Technology, PO Box 1334,
Durban 4000, Republic of South Africa}

\email{anpaliat@phys.uoa.gr}
\thanks{AP acknowledges the financial support of FONDECYT grant no. 3160121. }

\author{PGL Leach}
\address{School of Mathematical Sciences, University of KwaZulu-Natal, Durban, South Africa and}
\curraddr{Institute of Systems Science, Durban University of Technology, Durban, South Africa}
\email{leachp@ukzn.ac.za}
\thanks{PGLL acknowledges the support of the National Research Foundation of South Africa, the University of KwaZulu-Natal and the Durban University of Technology and thanks the Department of Mathematics, Pondicherry University, for gracious hospitality.}



\date{04:07:2018}


\begin{abstract}
The symmetry analysis of the Cheng Equation is performed. The Cheng Equation is reduced to a first-order equation of either Abel's Equations, the analytic solution of which is given in terms of special functions. Moreover, for a particular symmetry the system is reduced to the Riccati Equation or to the linear nonhomogeneous equation of Euler type. Henceforth, the  general solution  of the Cheng Equation with the use of the Lie theory is discussed, as also the application of Lie symmetries in a generalized Cheng equation.

\end{abstract}

\maketitle

\vspace{1.5cc}

\section{Introduction}
In 1984 Cheng\cite{cheng 03} discussed a pair of nonlinear partial differential equations with applications in the field of photosensitive molecules, namely
\begin{eqnarray}\label{1.1}
u_x & = & -auv\nonumber\\
v_t & = & bu_x ,
\end{eqnarray}
where $ u(t,x)$  and $v(t,x)$ represent the light intensity and density of the molecules, respectively.  More specifically, the first equation says that the amount of light absorbed by molecules in a film irradiated by a light beam
is proportional to the product of the light  intensity and density of the
molecules \cite{cheng 03}, while the second equation reads that the density
of the molecules in time is proportional to the light absorbed. The
constants $a$ and $b$ describe the absorption and proportionality constants,
respectively. The importance of the latter system is that it provides the first
example of exact solitary waves solutions on a failing film, where the kind of
waves are different from the solitons, observed experimentally by
Kapitza \& Kapitza in 1949 \cite{kap1}. Some later studies on the subject
are \cite{rev2,sw2,sw1}; for a review on the solitary waves we refer the
reader in \cite{rev1}. While a detailed review on solitary waves on failing
film is given in \cite{rev2}.

In the past work had been done to determine the solution of the Cheng equation through Hirota's bilinearisation\cite{cheng 01}  method and Painlev\'e analysis\cite{cheng 02}. The general travelling wave solution of (\ref{1.1}) had also been discussed in \cite{cheng 02}. In this paper we revisit this equation using Lie's method of symmetry analysis and discuss the possible reductions of the system and their solutions. The underlying Lie algebraic structure is also discussed. The Cheng Equation reduces to Abel's equations of the First and Second kinds, the solution of which can be given in terms of special functions and to the Riccati Equation and to a linear nonhomogeneous equation of Euler type.  To the authors' knowledge the reduction of the Cheng equation to Abel's  and linear Euler-type equation  have not been discussed in the literature and hence forms our main result. It is worthwhile mentioning that we found a new form of general solution which cannot follow from the general solutions as discussed before in the literature. Lie's approach makes the analysis more complete. The paper also discusses the case where the parameters of the Cheng Equation are space dependent.

The paper is arranged as follows.  In Section $2$ the Lie symmetries of (\ref{1.1})  are given. In subsequent sections the symmetry analysis corresponding to various cases of the arbitrary functions in the point symmetries is  discussed along with the case for which the parameters are space dependent.   The Conclusion and proper references are mentioned subsequently.

\section{Symmetry Calculation of the Cheng Equation }
For the convenience of the reader, we give a briefly discussion in the
theory of Lie point symmetries. In particular, we present the basic
definitions and main steps for the determination of Lie point symmetries for
a give differential equation.

Consider $H^{A}\left( t,x,u^{A},u_{,i}^{A}\right) =0$, to be a set of
differential equations, where $u^{A}=\left( u,v\right) ~$and $u_{,i}^{A}=%
\frac{\partial u^{A}}{\partial y^{i}}$ in which $y^{i}=\left( t,x\right) $.
Then under the action of the infinitesimal one-parameter point
transformation
\begin{eqnarray}
t^{\prime } &=&t\left( t,x,u^{A};\varepsilon \right)   \label{ls.01} \\
x^{\prime } &=&x\left( t,x,u^{A};\varepsilon \right)   \label{ls.02} \\
u^{^{\prime }A} &=&u^{A}\left( t,x,u^{A};\varepsilon \right)   \label{ls.03}
\end{eqnarray}%
in which $\varepsilon $ is an infinitesimal parameter, the set of
differential equations $H^{A}$ is invariant iff
\begin{equation}
H^{A}\left( t^{\prime },x^{\prime },u^{\prime A}\right) =H^{A}\left(
t,x,u^{A}\right)   \label{ls.04}
\end{equation}%
or equivalently \cite{bluman1}%
\begin{equation}
\lim_{\varepsilon \rightarrow 0}\frac{H^{A}\left( t^{\prime },x^{\prime
},u^{\prime A};\varepsilon \right) -H^{A}\left( t,x,u^{A}\right) }{%
\varepsilon }=0.  \label{ls.05}
\end{equation}%
The later expression is the definition of the Lie derivative $\mathcal{L}~$%
of $H^{A}$ along the direction
\begin{equation}
\Gamma =\frac{\partial t^{\prime }}{\partial \varepsilon }\partial _{t}+%
\frac{\partial x^{\prime }}{\partial \varepsilon }\partial _{x}+\frac{%
\partial u^{A}}{\partial \varepsilon }\partial _{u^{A}}.  \label{ls.06}
\end{equation}

Hence, we shall say that the vector field $\Gamma $ will be a Lie point
symmetry for the set of differential equations $H^{A}$ if and only if the
following condition is true%
\begin{equation}
\mathcal{L}_{\Gamma }\left( H^{A}\right) =0.  \label{ls.07}
\end{equation}

From (\ref{ls.07}) a set of linear differential equations is given for the
functions $\xi ^{t}=\frac{\partial t^{\prime }}{\partial \varepsilon },~\xi
^{x}=\frac{\partial x^{\prime }}{\partial \varepsilon }~$\ and $\eta ^{A}%
\frac{\partial u^{A}}{\partial \varepsilon }$, whose solution determine the
explicit form of the Lie point symmetries.

We omit the calculations and we give that the application of the Lie
symmetry condition (\ref{ls.07}) for the system (\ref{1.1}) provides with
the Lie symmetries
\begin{eqnarray}\label{ls.08}
\Gamma _{1} &=&-vg^{\prime }(x)\partial _{v}+g(x)\partial _{x}  \nonumber \\
\Gamma _{2} &=&h(t)\partial _{t}-uh^{\prime }(t)\partial _{u},  \nonumber\\
\end{eqnarray}
where $h(t)$ and $g(x)$ are arbitrary functions. The admitted Lie algebra of
system (\ref{1.1}) is $2A_{1}$ \footnote{%
We use the Mubarakzyanov Classification Scheme \cite{cheng 09,cheng 10,cheng
11,cheng 12}.}.\newline
We consider various possibilities for the arbitrary functions and conduct the reductions.
\begin{itemize}
\item[{\bf a.}] h(t) and g(x) are constant functions, which provide the translation symmetries.
\item[{\bf b.}] h(t)  and g(x)  are the identity function, which provide the scaling symmetries.
\item[{\bf c.}] The general case.
\end{itemize}
In the next three sections the reduction of (\ref{1.1}) with respect to various cases is discussed.

\section{Case I: $h(t)$ and $g(x)$ are constant functions }

\begin{eqnarray}\label{3.1}
\Gamma_{1A} & = & \partial_x\nonumber\\
\Gamma_{2A} & = & \partial_t.\nonumber\\
\end{eqnarray}

 $c\Gamma_{1A}$ + $\Gamma_{2A}$ gives a travelling-wave solution.

The similarity variables are
\begin{eqnarray}\label{3.2}
f     &=& x-ct\nonumber\\
u(x,t)&=&w(f)\nonumber\\
v(x,t)&=&k(f),\nonumber\\
\end{eqnarray}
 where $f$ is the new independent variable and $w(f)$ and $k(f)$ are the new dependent variables.
 This  reduces system (\ref{1.1}) to
 \begin{eqnarray}\label{3.3}
 w'(f) &=& -a k(f) w(f)\nonumber\\
 ck'(f) &=& -bw'(f).\nonumber\\
 \end{eqnarray}
System (\ref{3.3}) can be written as a second-order equation with respect to $w$, namely
\begin{equation}\label{3.4}
\frac{-c w'(f)^2}{a w(f)^2} + \frac{c w''(f)}{a w(f)} = b w'(f).\\
\end{equation}
The symmetries of (\ref{3.4}) are
\begin{eqnarray}\label{3.5}
\Gamma_{1B}&=&\partial_f\nonumber\\
\Gamma_{2B}&=&f\partial_f-w\partial_w.\nonumber\\
\end{eqnarray}

We consider $\Gamma_{1B}$ for reduction. The canonical coordinates are $$ n=w(f) \quad\mbox{\rm and}\quad m(n)=\frac{1}{w'(f)}.$$ This reduces (\ref{3.4}) to
\begin{equation}\label{3.6}
m'(n)=-\frac{n m(n)^2 ab}{c}-\frac{m(n)}{n}.\\
\end{equation}
The differential invariants,
$$ n=w(f) \quad\mbox{\rm and}\quad m=w'(f),$$ reduce (\ref{3.4}) to
\begin{equation}\label{3.7}
m'(n)=\frac{nab}{c}+\frac{m(n)}{n}.\\
\end{equation}
Equation (\ref{3.6}) is a Riccati equation and equation (\ref{3.7}) is a linear equation of Euler type.
Next we consider $ \Gamma_{2B}$ for reduction. The canonical coordinates are $$ n= f w(f) \quad\mbox{\rm and}\quad m(n)=\frac{1}{f(f w'(f)+w(f))}.$$
These reduce (\ref{3.4}) to
\begin{equation}\label{3.8}
m'(n)=\frac{(n^3ab+cn^2)m(n)^3}{cn}+\frac{(-n^2ab-cn)m(n)^2}{cn}-\frac{m(n)}{n}.\\
\end{equation}
The differential invariants,
$$ n= f w(f) \quad\mbox{\rm and}\quad m(n)= f^2 w'(f)),$$ reduce (\ref{3.4}) to
\begin{equation}\label{3.9}
m'(n)=\frac{m(n)(n^2 ab+m(n)+2cn)}{cn (m(n)+n)}.\\
\end{equation}
Equations (\ref{3.8}) and (\ref{3.9}) are Abel's equations of first and second kinds, respectively.
The solution of (\ref{3.6}) is
\begin{eqnarray}\label{3.10}
m(n)&=&\frac{c}{(n ab+cC_0)n},\nonumber\\
\end{eqnarray}
where $ C_0 $  is an  arbitrary constant.
The solutions of (\ref{3.8}) and (\ref{3.9}) are given in terms of Lambert W functions.

We use (\ref{3.10}) to derive the solutions for system (\ref{1.1}).  They are
\begin{eqnarray}\label{3.11}
u(x,t)&=&-\frac{c C_0}{-e^{-(c C_0((x-t)+C_1)}+ab}\nonumber\\
v(x,t)&=&\frac{c C_0}{a(-1+ab e^{c C_0((x-t)+C_1)})}.\nonumber\\
\end{eqnarray}

\section{Case II: $h(t)$  and $g(x)$  are the identity function}

The Lie point symmetries are now
\begin{eqnarray}\label{4.1}
\Gamma_{1C} & = & x\partial_x-v\partial_v\nonumber\\
\Gamma_{2C} & = & t\partial_t-u\partial_u.\nonumber\\
\end{eqnarray}

The similarity variables are
\begin{eqnarray}\label{4.2}
f     &=& \frac{t}{x},\nonumber\\
u(x,t)&=&\frac{w(f)}{t}\quad\mbox{\rm and}\nonumber\\
v(x,t)&=&\frac{k(f)}{x},\nonumber\\
\end{eqnarray}
 where $f$ is the new independent variable and $w(f)$ and $k(f)$ are the new dependent variables.
 These  reduce system (\ref{1.1}) to
 \begin{eqnarray}\label{4.3}
 fw'(f) &=& a k(f) w(f)\nonumber\\
 k'(f) &=& -bw'(f).\nonumber\\
 \end{eqnarray}
Equations (\ref{4.3}) can be written as a second-order equation with respect to $w$.  It is
\begin{equation}\label{4.4}
\frac{w'(f)}{a w(f)}+ \frac{f w''(f)}{a w(f)} -\frac{f w'(f)^2}{a w(f)^2} + b w'(f) = 0.\\
\end{equation}
The symmetries of (\ref{4.4}) are
\begin{eqnarray}\label{4.5}
\Gamma_{1D}&=&\partial_f\nonumber\\
\Gamma_{2D}&=&-f\log(f)\partial_f+w\partial_w.\nonumber\\
\end{eqnarray}

The analysis is similar to the previous case and the equation reduces to Riccati, Linear Euler  and Abel's of the first and second kind, respectively.
Next we perform the analysis by another set of similarity variables.

They are
\begin{eqnarray}\label{4.6}
f     &=& \frac{t}{x},\nonumber\\
u(x,t)&=&\frac{w(f)}{x}\quad\mbox{\rm and}\nonumber\\
v(x,t)&=&\frac{k(f)}{x},\nonumber\\
\end{eqnarray}
 where $f$ is the new independent variable and w(f) and k(f) are the new dependent variables.
 These  reduce system (\ref{1.1}) to
 \begin{eqnarray}\label{4.7}
 f w'(f)+w(f) &=& a k(f) w(f)\nonumber\\
 k'(f) &=& -b(f w'(f)+w(f)).\nonumber\\
 \end{eqnarray}
The system (\ref{4.7}) can be written as a second-order equation in $w$, namely
\begin{equation}\label{4.8}
f w(f) w''(f)+ a b w(f)^3 + w(f) w'(f) +a b f w(f)^2 w'(f) - f w'(f)^2 = 0.\\
\end{equation}
The symmetries of (\ref{4.8}) are
\begin{eqnarray}\label{4.9}
\Gamma_{1E}&=&f\partial_f-w\partial_w\quad\mbox{\rm and}\nonumber\\
\Gamma_{2E}&=&f\log(f)\partial_f-(\log(f)+1)w\partial_w.\nonumber\\
\end{eqnarray}

 The canonical coordinates with respect to $\Gamma_{1E}$ are
 \begin{eqnarray}\label{4.10}
 r&=&w(f)\quad\mbox{\rm and}\nonumber\\
 v(r)&=&\frac{1}{f(fw'(f)+w(f))}.\nonumber\\
 \end{eqnarray}

Equation (\ref{4.8}) reduces to
 \begin{equation}\label{4.11}
 v'(r)= v(r)^2 r a b-\frac{v(r)}{r}.
 \end{equation}
 The differential invariants with respect to $\Gamma_{1E}$ are
 \begin{eqnarray}\label{4.12}
 r&=&fw(f)\quad\mbox{\rm and}\nonumber\\
 v(r)&=&f^2 w'(f).\nonumber\\
 \end{eqnarray}
Equation (\ref{4.8}) reduces to
 \begin{equation}\label{4.13}
 v'(r)=-rab+\frac{v(r)}{r}.
 \end{equation}
Equation (\ref{4.11}) is a Riccati equation and equation (\ref{4.13}) is an equation of  Euler type.
 The canonical coordinates with respect to $\Gamma_{2E}$ are
 \begin{eqnarray}\label{4.14}
 r&=&w(f)f\log(f)\quad\mbox{\rm and}\nonumber\\
 v(r)&=&\frac{1}r\nonumber
 \end{eqnarray}

Equation (\ref{4.8}) reduces to
  \begin{equation}\label{4.15}
  v'(r) = -\frac{(abr^3-r^2)v(r)^3}{r}-\frac{(-abr^2+r)v(r)^2}{r}-\frac{v(r)}{r}.
  \end{equation}
 The differential invariants with respect to $\Gamma_{2E}$ are
 \begin{eqnarray}\label{4.16}
  r&=&w(f)f\log(f)\quad\mbox{\rm and}\nonumber\\
 v(r)&=&(w(f)\log(f)+w(f)+f\log(f)w'(f))f\log(f).\nonumber\\
\end{eqnarray}
Equation (\ref{4.8}) reduces to
\begin{equation}\label{4.17}
v'(r)=\frac{v(r)}{r}-rab+1+\frac{abr^2-r}{v(r)}.
\end{equation}
Equations (\ref{4.14}) and (\ref{4.16}) are Abel's equation of First kind and Second kind, respectively.
the solution of which can be given in terms of Lambert W function.

\section{The general case}

The general symmetry is considered.

The similarity variables are
\begin{eqnarray}\label{5.1}
f     &=& h_{a}(t)-g_{a}(x),\nonumber\\
u(x,t)&=&h_{a}'(t)w(f)\quad\mbox{\rm and}\nonumber\\
v(x,t)&=&g_{a}'(x)k(f),\nonumber\\
\end{eqnarray}
 where $f$ is the new independent variable, $w(f)$ and $k(f)$ are the new dependent variables and  $h_{a}(t)$ and $g_{a}(x)$  are given as
 \begin{eqnarray}\label{5.2}
 h_{a}(t)&=&\int \frac{1}{h(t)} dt\quad\mbox{\rm and}\nonumber\\
 g_{a}(x)&=&\int \frac{1}{g(x)} dx.\nonumber\\
 \end{eqnarray}
 These reduce (\ref{1.1}) to
 \begin{eqnarray}\label{5.3}
 w'(f)&=&aw(f)k(f)\quad\mbox{\rm and}\nonumber\\
 k'(f)&=&-bw'(f).\nonumber\\
 \end{eqnarray}
 This is similar to the analysis of Case I. When the similarity variables are chosen as in (\ref{4.2}), the equation reduces to Abel's Equation of the  First and Second kind,  Riccati and nonhomogeneous equation of Euler type. The general solution is given as
 \begin{eqnarray}\label{5.4}
u(x,t)&=&\frac{ C_0 h_{a}'(t)}{e^{- C_0((x-t)+C_1)}+ab}\quad\mbox{\rm and}\nonumber\\
v(x,t)&=&-\frac{ C_0 g_{a}'(x)}{a(-1+ab e^{ C_0((x-t)+C_1)})}.\nonumber\\
\end{eqnarray}
  Similarly, when the similarity variables are chosen as in (\ref{4.6}), i.e.,
 \begin{eqnarray}\label{5.5}
f     &=& h_{a}(t)-g_{a}(x),\nonumber\\
u(x,t)&=&g_{b}(x)w(f)\quad\mbox{\rm and}\nonumber\\
v(x,t)&=&g_{a}'(x)k(f),\nonumber\\
\end{eqnarray}
 where $h_{a}(t)$ and $g_{a}(x)$ are as defined in (\ref{5.2}) and $g_{b}(x)$ is given in terms of $h(t)$ and $g(x).$ The results follows as in Case II.\\
\section{parameters are space dependent}
 From the first equation in equation (\ref{1.1}), we have
\[
-\frac{1}{a}\frac{u_{x}}{u}=v
\]%
so that%
\[
v_{t}=-\frac{1}{a}\left( \frac{u_{x}}{u}\right) _{t}=-\frac{1}{a}\left(
\frac{u_{xt}}{u}-\frac{u_{x}u_{t}}{u^{2}}\right).
\]

Finally%
\[
\left( \frac{u_{xt}}{u}-\frac{u_{x}u_{t}}{u^{2}}\right) +\frac{b}{a}%
u_{x}=0.
\]%
The symmetry vectors are
\[
\lambda \left( x\right) \partial _{x}~,~\tau \left( t\right) \partial
_{t}-\tau ^{\prime }\left( t\right) u\partial _{u}.
\]

\bigskip

If $a=a\left( x\right) $ and $b=b\left( x\right) $, then%
\begin{equation}\label{6.1}
\left( \frac{u_{xt}}{u}-\frac{u_{x}u_{t}}{u^{2}}\right) +c\left( x\right)
u_{x}=0
\end{equation}
where $c\left( x\right) =\frac{b}{a}$, and the symmetry vectors are
\[
c\left( x\right) \partial _{x}-c^{\prime }\left( x\right) u\partial
_{u}~,~\tau \left( t\right) \partial _{t}-\tau ^{\prime }\left( t\right)
u\partial _{u}.
\]
 The similarity variable with respect to the first symmetry vector is
 $u(x,t)= \frac{s(x)}{c(x)},$ where $s(x)$ is the new dependent variable.
 When $s(x)=x$, the equation(\ref{6.1})  reduces to the Euler equation,
 \begin{equation}\label{6.2}
 c'(x)-\frac{c(x)}{x}=0.
 \end{equation}
 In general, when $s(x)$ is a polynomial function, equation(\ref{6.1}) reduces to a linear Euler equation.

 For any $s(x)$, the equation(\ref{6.1}) reduces to
 \begin{equation}\label{6.3}
 \frac{c'(x)}{s'(x)}-\frac{c(x)}{s(x)}=0.
 \end{equation}
The linear combination of the two symmetry vectors is
$$
c\left( x\right) \partial _{x}-c^{\prime }\left( x\right) u\partial
_{u}~+C_{2}(~\tau \left( t\right) \partial _{t}-\tau ^{\prime }\left( t\right),
u\partial _{u})
$$
where $C_{2}$ is an arbitrary constant.
The similarity variables are
\begin{eqnarray}\label{6.4}
f_{a}&=& c_{a}(x)-\tau_{a}(t) \quad\mbox{\rm and}\nonumber\\
u(x,t)&=&w_{a}(f_{a}),\nonumber\\
\end{eqnarray}
 where $f_{a}$ is the new independent variable, $w_{a}(f_{a})$ is the new dependent variables and  $c_{a}(x)$ and $\tau_{a}(t)$  are given as
 \begin{eqnarray}\label{6.5}
 c_{a}(x)&=&\int \frac{1}{c(x)} dx\quad\mbox{\rm and}\nonumber\\
 \tau_{a}(t)&=&\int \frac{1}{\tau(t)} dt.\nonumber\\
 \end{eqnarray}
For $ c_{a}(x)=x$ and $h_{a}(t)=t$, the equation(\ref{6.1}) reduces to
\begin{equation}\label{6.6}
\frac{w_{a}''(f_{a})}{w_{a}(f_{a})}=\frac{w_{a}'(f_{a})}{w_{a}(f_{a})^2}+w_{a}'(f_{a}).
\end{equation}
The Lie symmetries are
\begin{eqnarray}\label{6.7}
\Gamma_{a}&=&\partial_{f_{a}}\nonumber\\
\Gamma_{b}&=&f_{a}\partial_{f_{a}}-w_{a}\partial_{w_{a}}.
\end{eqnarray}
The equation (\ref{6.6}) behaves similarly to(\ref{3.4}).  The symmetry $\Gamma_{a}$ reduces (\ref{6.6}) to the Riccati equation and $\Gamma_{b}$ reduces(\ref{6.6}) to Abel's Equation.

\section{conclusion}
In this work we applied the theory of invariant transformations to the study and the determination of analytical solutions for the Cheng equation. We found that the similarity solutions are expressed in terms of solutions of the Abel equations, that is, Cheng equation, when it is reduced from a second-order partial differential equation to a first-order ordinary differential equation with the use of point transformations, is equivalent with the Abel equation.

However, we found that the symmetry vectors depend upon arbitrary functions and we have performed our analysis for  for  various choices of similarity variable.   The results of the various cases are similar, as for all the cases the system is reduced to either an  Abel's or a Riccati equation. The paper also made a note of the Lie algebra.  Finally, the case when the parameters are space dependent is also discussed. The general solution of the later equation have been presented by the Lie symmetry approach.

\end{document}